\documentclass{article}
\usepackage{a4}
\usepackage[latin1]{inputenc}
\usepackage[T1]{fontenc}
\usepackage{amsmath}
\usepackage{amssymb}
\usepackage{amsthm}
\usepackage[ps,dvips,arrow,matrix,tips]{xy}
\SelectTips{cm}{10}

\renewcommand{\phi}{\varphi}
\renewcommand{\emptyset}{\varnothing}
\newtheorem{thm}{Theorem}[section]
\newtheorem{prop}[thm]{Proposition}
\newtheorem{lemma}[thm]{Lemma}
\newtheorem{defn}[thm]{Definition}
\newenvironment{pf}[1][]{\noindent{\emph{Proof}#1.} }{\hfill $\square$}
\newenvironment{sketchofproof}[1][]{\noindent{\emph{Sketch of proof}#1.} }{\hfill $\square$}

\DeclareMathOperator{\Spec}{Spec}
\renewcommand{\P}{{\mathbf{P}}}

\newcommand{\Aut}{\mathrm{Aut}}
\newcommand{\F}[1]{\mathbf{F}_{#1}}
\newcommand{\Falg}{\mathbf{F}}
\newcommand{\Fq}{{\F{q}}}
\newcommand{\N}{{\mathbf{N}}}
\newcommand{\Z}{{\mathbf{Z}}}
\newcommand{\C}{{\mathbf{C}}}
\newcommand{\Gal}{\mathrm{Gal}}
\newcommand{\tuple}[2]{#1, \mskip2.5mu \ldots \mskip-1mu, \mskip2.5mu #2}
\newcommand{\symm}[1]{\mathfrak{S}(#1)}
\newcommand{\abs}[1]{\left| \mskip1mu #1 \right|}

\title{A potential analogue of Schinzel's hypothesis for polynomials with coefficients in $\Fq[t]$}
\date{July 1, 2005}
\author{Andreas O. Bender\thanks{The first author acknowledges financial support provided by the Japan Society for the Promotion of Science (JSPS) as well as through the European Community's Human Potential Programme under contracts HPRN-CT-2000-00120 [AAG] and HPRN-CT-2000-00114 [GTEM].}\\
Korea Institute for Advanced Study\\
Seoul 130-722\\
Korea\\
\texttt{andreas@kias.re.kr} \and
Olivier Wittenberg\\
Laboratoire de mathématiques, Bâtiment 425\\
Université Paris-Sud\\
F-91405 Orsay\\
France\\
\texttt{olivier.wittenberg@ens.fr}}

\begin{document}
\maketitle

\section{Introduction}

The Schinzel hypothesis essentially claims that finitely many irreducible polynomials in one variable over $\Z$ simultaneously assume infinitely many prime values unless there is an obvious reason why this is impossible.

We prove that under a restriction on the characteristic and a smoothness assumption, finitely many irreducible polynomials in one variable over the ring $\Fq[t]$ assume simultaneous prime values after a sufficiently large
extension of the field of constants.

\subsection{The Schinzel hypothesis over $\Z$}
Let $\tuple{f_{1}(x)}{f_{r}(x)}$ be irreducible
polynomials with coefficients in~$\Z$. Assume that
the leading coefficient of every $f_{i}(x)$ is positive and that for each prime~$p$, there exists an integer $x_{p}$ such that no 
$f_{i}(x_{p})$ is divisible by $p$.
Then $\tuple{f_{1}(x)}{f_{r}(x)}$ are simultaneously prime for infinitely many integer values of $x$.\\
In its present generality, this conjecture was first stated in~\cite{schinzel}.

\subsection{The Schinzel hypothesis over $\Fq[t]$}
A naïve analogon to the Schinzel hypothesis over the coefficient ring $\Fq[t]$ can be formulated as follows:

\vspace{1mm}
Let $\tuple{f_{1}(x)}{f_{r}(x)}$ be non-constant polynomials in~$\Fq[t,x]$ which are irreducible in~$\Fq(t)[x]$ and assume that 
for each prime ${\mathfrak p}$ of $\Fq[t],$ there is an $x_{\mathfrak p}$ in~$\Fq[t]$  
such that no $f_{i}(x_{\mathfrak p})$ is divisible by ${\mathfrak p}$.
Then $\tuple{f_{1}(x)}{f_{r}(x)}$ are simultaneously prime for infinitely many values of $x$ in $\Fq[t]$.

\vspace{1mm}
In this form, the hypothesis is known to be false; counterexamples
are described in~\cite{CConrad}, one of them being the following: With $p$ the characteristic of $\Fq$, choose an integer $b$ with $1<b<4q$ and $(b,p(q-1))=1$ (e.g., $b=2q-1$). Let $f(x)=x^{4q}+t^{b}$. Then 
$f(g)$ is reducible for all $g\in \Fq[t]$ (see~\cite[Section~4]{CConrad}).

In~\cite{bateman}, Bateman and Horn formulated a quantitative version of the Schinzel hypothesis, which specifies the conjectural proportion of integers for which the polynomials assume prime values.

In~\cite[(1.2),(1.8)]{CConrad}, Conrad, Conrad and Gross presented an analogue to the Bateman--Horn conjecture for the case of one polynomial over $\Fq[t]$, supported by good agreement with numerical evidence. As for the qualitative case, this conjecture implies that the naïve function field version of the Schinzel hypothesis enunciated above does in fact hold for one separable polynomial.

For polynomials with coefficients in either $\Z$ or $\Fq[t]$, Dirichlet's theorem about primes in arithmetic progressions and its analogue for function fields is the only case in which the Schinzel hypothesis is known to hold; see~\cite{kornblum} or~\cite[Theorem~4.7]{rosen} for the case of $\Fq[t]$. Note that in that case of one polynomial of degree 1, these results amount to the quantitative statements of the conjectures of Bateman--Horn and of Conrad, Conrad and Gross, respectively.

The goal of this note is to prove the following theorem.

\begin{thm}
\label{mainthm}
Let $\Fq$ be a finite field of characteristic $p$ and cardinality $q$.
Let $\tuple{f_1}{f_n} \in \Fq[t,x]$ be irreducible polynomials whose
total degrees $\deg(f_i)$ satisfy $p \nmid \deg(f_i) (\deg(f_i)-1)$
for all~$i$. Assume that the curves $C_i \subset \P^2_\Fq$ defined
as the Zariski closures of the affine curves
$$
f_i(t,x)=0
$$
are smooth. Then, for any sufficiently
large $s \in \N$, there exist \mbox{$a, b \in \F{q^s}$} such that
the polynomials $\tuple{f_1(t,at+b)}{f_n(t,at+b)} \in \F{q^s}[t]$
are all irreducible.
\end{thm}

\vspace{3mm}
{\bf Acknowledgements}: 
A preliminary manuscript about this topic, written solely by the first author,
contained some gaps; he is grateful to Jean-Louis Colliot-Thélène, Jürgen
Klüners and Sir Peter Swinnerton-Dyer for pointing them out.

We had very useful discussions about the contents of this paper with Jean-Louis Colliot-Thélène; previous versions benefited from comments made by Ido Efrat, Bert van Geemen, Pierre Dèbes, Moshe
Jarden, and Fumiharu Kato.
We also thank Keith Conrad for corresponding with us about the results contained in~\cite{CConrad}
and Ofer Gabber for pointing out that
Lemma~\ref{serregeninertia} could be used instead of~\cite[Proposition~4.4.6]{serre}.

This research was carried out while the first author was staying at the University of Pavia and Collegio Ghislieri, the University of Padova, at Ben-Gurion University of the Negev, 
the Hebrew University of Jerusalem, with the major part having been done at Kyoto University. To all these institutions, he expresses his gratitude for their hospitality.

\vspace{3mm}
\noindent\emph{Notations.}
We denote by $\abs{S}$ the cardinality of a set $S$ and by $\symm{S}$ the
symmetric group on $S$. Let $k$ be a field and $X$ be a $k$-scheme. If $k'$
is a field extension of $k$, the scheme $X \times_{\Spec(k)} \Spec(k')$ will often be denoted $X_{k'}$.
We write $\kappa(x)$ for the residue field of $x \in X$, and $\kappa(X)$ for the
function field of $X$ when $X$ is integral. Finally, when $Y$ is an $X$-scheme,
$\Aut_X(Y)$ denotes the group of $X$-automorphisms of~$Y$.

\section{A preliminary result about generic covers of $\P^1$}

If $k$ is a field, a finite $k$-scheme $X$ will be said to have \emph{at most one double point} if
$n(X) \geq r(X)-1$, where $r(X)$ and $n(X)$ respectively denote
the rank and the geometric number of points of~$X$.

The following definition was introduced by Hurwitz~\cite{hurwitz} in his proof of connectedness of the moduli spaces for curves of genus $g$ over $\C$. 

\begin{defn}
A finite morphism $f\colon C \rightarrow \P^{1}_{k}$ is called generic if $f^{-1}(x)$ has at most one double point for all $x\in \P^{1}_{k}$.
\end{defn}

Note that if the characteristic of~$k$ is not~$2$, $f$ is generic and~$C$
is an integral curve, then~$f$ is separable (in the sense that the field
extension $\kappa(C)/\kappa(\P^1_k)$ is separable).

\begin{prop}
\label{prelprop}
Let $C$ be a regular, complete, geometrically irreducible curve over a field~$k$,
endowed with a finite separable generic morphism $f \colon C \rightarrow \P^1_k$.
Let $C'$ be a regular, complete, irreducible curve over $k$, and $g \colon C' \rightarrow C$
be a finite morphism. Assume that the finite extension $\kappa(C')/\kappa(\P^1_k)$ is a Galois
closure of the subextension $\kappa(C)/\kappa(\P^1_k)$.
We denote respectively by $G$ and $H$ the Galois groups of
$\kappa(C')/\kappa(\P^1_k)$ and $\kappa(C')/\kappa(C)$.
Then $C'$ is geometrically irreducible over~$k$ and the morphism
$$
G \longrightarrow \symm{H \backslash G}
$$
induced by right multiplication is an isomorphism.
Moreover, all the ramification indices of $\kappa(C')/\kappa(\P^1_k)$
are $\leq 2$.
\end{prop}

\begin{pf}
Let $k'$ denote the algebraic closure of $k$ in $\kappa(C')$.
We denote respectively by $G'$ and $H'$ the subgroups of $G$ defined
by the subfields $\kappa(\P^1_{k'})$ and $\kappa(C_{k'})$ of $\kappa(C')$, so that
we have a canonical commutative diagram as follows, where the labels indicate the Galois
groups of the generic fibres:
$$
\xymatrix@C=8ex{
C_{k'} \ar[dd] \ar[rr] && \P^1_{k'} \ar[dd] \\
& C' \ar[ul]_{H'} \ar[ur]^{G'} \ar[dl]_H \ar[dr]^G \\
C \ar[rr] && \P^1_k
}
$$
Let us endow $H \backslash G$ (resp.~$H' \backslash G'$) with the action of~$G$ (resp.~$G'$)
by right multiplication.
The equality $H \cap G' = H'$ of subgroups of $G$ yields a natural injective $G'$-equivariant map
$H' \backslash G' \rightarrow H \backslash G$, which is even bijective since
$\abs{H \backslash G}$ and $\abs{H' \backslash G'}$ are both equal to~$\deg(f)$.
Hence a commutative square
$$
\special{ps: TeXDict begin @defspecial /XYdict where pop begin XYdict begin
/fullhook{1 .8 scale 0 xysegl .5 mul 2 div dup -90 90 arcn 1 1 .8 div
scale 1 setlinecap}bind def end end @fedspecial end}%
\xymatrix{
G' \vphantom{\backslash} \ar[r] \ar@{^{(}->}[d] & \symm{H' \backslash G'} \ar@{=}[d] \\
G \vphantom{\backslash} \ar[r] & \symm{H \backslash G} \rlap{\text{,}}
}
$$
where the horizontal arrows are induced by right multiplication.
The bottom horizontal arrow is injective, in virtue of the equality
$\bigcap_{a \in G} aHa^{-1} = 1$, itself a consequence of the
hypothesis that
$\kappa(C')/\kappa(\P^1_k)$ is a Galois closure of $\kappa(C)/\kappa(\P^1_k)$.
For the first part of the proposition, it
only remains to be shown that the top horizontal arrow is surjective;
indeed, this will imply not only that the bottom horizontal arrow is
an isomorphism, but also that $G'=G$, hence $k'=k$, which is equivalent
to $C'$ being geometrically irreducible over $k$. The second part will follow from
the first once we know that all ramification indices of $\kappa(C')/\kappa(\P^1_{k'})$
are $\leq 2$.

We shall now make use of the following classical result.
A very similar lemma is stated and proven in \cite[Proposition~4.4.6]{serre}.

\begin{lemma}
\label{serregeninertia}
Let $X$ be a regular, complete, geometrically irreducible curve over a field~$K$, endowed
with a finite and generically Galois morphism $X \rightarrow \P^1_K$ with group~$G$.
Then~$G$ is generated by the inertia subgroups above closed
points of $\P^1_K$ and their conjugates.
\end{lemma}

\begin{proof}
Let $H \subseteq G$ denote the normal subgroup generated by the inertia subgroups
and their conjugates.
The map $X \rightarrow \P^1_K$ can be factored
as $X \rightarrow Y \rightarrow \P^1_K$, where $Y$ is a regular, complete, irreducible
curve over~$K$ whose function field is the subfield of~$\kappa(X)$ fixed by~$H$.
The curve~$Y$ is geometrically irreducible over~$K$, since~$X$ is, and it follows
from Lemma~\ref{wellknownlemma} that it is étale over~$\P^1_K$;
therefore $Y=\P^1_K$ (see \cite[IV.2.5.3]{hartshorne}), hence $H=G$.
\end{proof}

Let us consider the cover $C' \rightarrow \P^1_{k'}$. It is generically Galois with
group $G'$. Let $I \subseteq G'$ be the inertia subgroup of $G'$ associated with a point
of $C'$ whose image by $f \circ g$ will be denoted $x$.
By Lemma~\ref{wellknownlemma},
the geometric number of points of $f^{-1}(x)$ is $\abs{H' \backslash G'/I}$.
Moreover, the rank of $f^{-1}(x)$ is $\abs{H' \backslash G'}$. The hypothesis that
$f^{-1}(x)$ has at most one double point thus leads to the inequality
$$\abs{H' \backslash G'/I} \geq \abs{H' \backslash G'} - 1\text{,}$$
thereby proving that every non-trivial inertia subgroup of $G'$ has order~$2$ and acts
as a transposition
on $H' \backslash G'$.
Applying Lemma~\ref{serregeninertia} to $X=C'$ and $K=k'$ now yields that
$G'$ is generated by elements which act on $H' \backslash G'$ as transpositions. The image
of $G' \rightarrow \symm{H' \backslash G'}$ is therefore a transitive subgroup
of $\symm{H' \backslash G'}$ which is
generated by transpositions; but the only such subgroup is $\symm{H' \backslash G'}$
itself (see~\cite[Lemma~4.4.4]{serre}), hence the result.
\end{pf}

\section{Proof of Theorem~\ref{mainthm}}

To prove Theorem~\ref{mainthm}, we may and will
assume that the polynomials $(f_i)_{1 \leq i \leq n}$ are pairwise non-proportional.
Let~$\Falg$ denote an algebraic closure of~$\Fq$. The symbol~$\F{q^s}$ will now be understood
to refer to the unique subfield of $\Falg$ with cardinality $q^s$. Let~$M_0 \in \P^2(\Fq)$
denote the point at infinity with coordinates $x=1$, $t=0$.

\begin{prop}
\label{existsu}
There exists a non-empty open subset $U \subset \P^2_{\Fq} \setminus \{M_0\}$, disjoint from~$C_i$
for all $i \in \{\tuple{1}{n}\}$, such that every line~$D$ in $\P^2_\Falg$ which
meets~$U$ satisfies the following properties:
\begin{enumerate}
\item For all $i \in \{\tuple{1}{n}\}$, the scheme-theoretic intersection $(C_i)_\Falg \cap D$
has at most one double point (as a finite $\Falg$-scheme).
\item The line $D$ is not tangent to more than one of the curves $(C_i)_\Falg$, $i \in \{\tuple{1}{n}\}$.
\end{enumerate}
\end{prop}

\begin{pf}
It is enough to prove that there are finitely many lines~$D$ in $\P^2_\Falg$ not satisfying
the above properties. Indeed, once this is known, we can take for $U$ any non-empty
open subset disjoint from the curves $C_i$ and from all these lines.

We shall use the duality theory of plane curves.  To every irreducible
plane curve~$C \subset\P^2$ over some field is associated an
irreducible curve $C^\star \subset (\P^2)^\star$, called its dual,
together with a canonical rational map $\rho \colon C \dashrightarrow
C^\star$, called the Gauss map, which sends a smooth point of~$C$ to
its tangent line. Here $(\P^2)^\star$ denotes the dual projective
plane. The reader is referred to~\cite{kleimanvanc} for an overview of
this theory, of which we shall only use the following two results. Firstly, the
Monge-Segre-Wallace criterion (see \cite[p.~169]{kleimanvanc}) ensures that
the equality $C^{\star \star} = C$ of closed subsets of~$\P^2$ holds as soon
as $\rho$ is separable. Secondly, it follows from Corollaire~3.5.0 and Corollaire~3.2.1
of \cite{katz} that if $\rho$ is separable, then there are only finitely many lines~$D$ in $\P^2$ such
that the scheme-theoretic intersection $C \cap D$ does not have at most one double point.

Let us apply these results to the smooth curves $(C_i)_\Falg \subset
\P^2_\Falg$.  A line in $\P^2_\Falg$ which is tangent to more than one
of the curves $(C_i)_\Falg$, $i \in \{\tuple{1}{n}\}$, corresponds to
a point in $(C_i)_\Falg^\star \cap (C_j)_\Falg^\star$ for distinct $i,
j$. Now if the Gauss maps $\rho_i \colon C_i \rightarrow C_i^\star$
and $\rho_j \colon C_j \rightarrow C_j^\star$ are separable, this
intersection is finite. Indeed, the curves $C_i^\star$ and $C_j^\star$
being irreducible, they would otherwise be equal; but the separability
of $\rho_i$ and $\rho_j$ implies that $C_i^{\star\star}=C_i$ and
$C_j^{\star\star}=C_j$, and we have assumed that $C_i\neq C_j$.

We are thus reduced to proving that the maps $\rho_i$ are all
separable.  As the curves~$C_i$ are smooth, Euler's formula shows that the Gauss maps $\rho_i$ can be extended to morphisms $r_{i}:\P^{2}\rightarrow \P^{2}$. An application of the projection formula for intersections \cite[A.1, A4]{hartshorne} then gives
$$r_{i\star} (C_i.r_{i}^{\star} D)=(r_{i\star} C_i).D,$$
where $D$ is a line in the target space $\P^{2}$. The equations of $r_{i}$ now show that $r_{i}^{\star} D$ has degree $\deg(C_i)-1$. By definition of $r_{i\star}$, we have $r_{i\star} C_i=\deg(\rho_i) C_i^{\star}$ and so Bézout's theorem implies the following formulae:
\begin{equation*}
\deg(C_i) \left(\deg(C_i)-1\right) = \deg(\rho_i) \deg(C^\star_i) \text{.}
\end{equation*}
The hypothesis on the total degrees of the polynomials $f_i$ now implies that $\deg(\rho_i)$ is
prime to~$p$, hence $\rho_i$ is separable.
\end{pf}

\bigskip
Let $U \subset \P^2_{\Fq}$ be given by Proposition~\ref{existsu} and $s_0 \in \N$ be large
enough so that $U(\F{q^s})\neq \emptyset$ for all $s \geq s_0$. Let $s \in \N$ be a sufficiently
large integer; for the time being, this means that $s \geq s_0$, but another condition on~$s$
will be introduced later.
For the sake of clarity,
we will henceforth denote the field $\F{q^s}$ by~$k$.
Fix $M \in U(k)$ and
denote by $\phi_i \colon (C_i)_k \rightarrow \P^1_k$ the
$k$-morphism obtained by composing the inclusion $(C_i)_k \subset
\P^2_k \setminus \{M\}$ with the morphism
$\P^2_k \setminus \{M\} \rightarrow \P^1_k$
defined by projection from~$M$.
The morphism $\phi_i$ is finite of degree $\deg(f_i)$ and is
generic, since $M \in U$. Being generic, it is separable
(note that the hypotheses of Theorem~\ref{mainthm} imply that~$p\neq 2$);
therefore there exists a smooth, complete,
connected curve $C'_i$ over~$k$ and a finite morphism $C'_i \rightarrow (C_i)_k$, such that
the induced field extension $\kappa(C'_i)/\kappa(\P^1_k)$ is a Galois closure of
$\kappa((C_i)_k)/\kappa(\P^1_k)$.
Let us write, for simplicity, $K=\kappa(\P^1_k)$,
$K_i=\kappa(C'_i)$, $G_i=\Gal(K_i/K)$ and $H_i=\Gal(K_i/\kappa((C_i)_k))$.
Proposition~\ref{prelprop} now shows that for all $i \in \{\tuple{1}{n}\}$, the curve $C'_i$ is
geometrically connected over~$k$,
the group $G_i$ is canonically isomorphic to $\symm{H_i \backslash G_i}$ and the ramification indices
of $K_i/K$ are $\leq 2$.

Let $R_i \subset \P^1_k$ denote the branch locus of the morphism $C'_i \rightarrow \P^1_k$.

\begin{prop}
\label{ridisjoint}
The subsets $R_i \subset \P^1_k$ for $i \in \{\tuple{1}{n}\}$ are pairwise disjoint.
\end{prop}

\begin{pf}
We shall need the following well-known lemma, which is a direct consequence
of Lemma~\ref{wellknownlemma}.

\begin{lemma}
Let $E/K$ be a finite separable extension of global fields, and let $L$ be a Galois closure
of $E/K$. Then a finite place of $K$ is unramified in $E$ if and only if it is unramified
in $L$.
\end{lemma}

The lemma shows that $R_i$ is also the branch locus of the morphism
$(C_i)_k \rightarrow \P^1_k$. An $\Falg$-point of $R_i \cap R_j$ therefore gives rise
to a line in $\P^2_\Falg$ which is both tangent to $(C_i)_\Falg$ and $(C_j)_\Falg$,
and which contains $M$.
As $M \in U$, there is no such line if $i \neq j$, hence the proposition.
\end{pf}

\bigskip
Let $L$ denote the ring $K_1 \otimes_K \cdots \otimes_K K_n$.

\begin{prop}
The ring $L$ is a field, and $k$ is separably closed in~$L$.
\end{prop}

\begin{pf}
For $j \in \{\tuple{0}{n}\}$, let us write $L_j$ for the $K$-algebra
$K_1 \otimes_K \cdots \otimes_K K_j$ and
prove that~$L_j$ is a field in which $k$ is separably closed,
by induction on $j$.
The case $j=0$ is trivial, as $L_0 = K$.
Assume now that $j>0$ and that~$L_{j-1}$ is a field in which~$k$ is separably closed.
Let $\Omega$ be a field containing~$\Falg$, $L_{j-1}$ and $K_j$.
Consider the subfield $E_j \subset \Omega$ defined as the intersection of
the composita $\Falg L_{j-1}$ and $\Falg K_j$.
Being a finite extension of $\Falg K$, it is the function field of a connected finite
cover of $\P^1_\Falg$. Proposition~\ref{ridisjoint} now shows that this cover is unramified; a
connected finite étale cover of $\P^1_\Falg$ is necessarily trivial, hence $E_j=\Falg K$.
As~$\Falg K_j$ is Galois and~$\Falg L_{j-1}$ is finite over $\Falg K$, this
is enough to imply that~$\Falg L_{j-1}$ and~$\Falg K_j$ are linearly disjoint subfields
of~$\Omega$ over $\Falg K$; in other words, $\Falg L_{j-1} \otimes_{\Falg K} \Falg K_j$ is
a field. We have~$\Falg L_{j-1} = \Falg \otimes_k L_{j-1}$ and $\Falg K_j = \Falg
\otimes_k K_j$ since~$k$ is separably closed in~$L_{j-1}$ and in $K_j$, hence
$\Falg L_{j-1} \otimes_{\Falg K} \Falg K_j = \Falg \otimes_k L_j$. As this ring is a field,
$k$ is separably closed in~$L_j$.
\end{pf}

\bigskip
Let $C'$ denote a smooth complete connected curve over~$k$ with function field $L$.
There is a natural finite morphism $\psi \colon C' \rightarrow \P^1_k$, which is
generically Galois and therefore separable.
We denote by $g$ the genus of $C'$, by $G$ the group $\Gal(L/K)$, by~$N$ the
degree of $\psi$, and by $(x, L/K)$ the Artin symbol of the extension $L/K$ above
a closed point $x \in \P^1_k$ which does not ramify in $L$.
We would now like to find a rational point of $\P^1_k$ above which the fibre of $\psi$ is
integral. To this end, we resort to an effective version of the \v{C}ebotarev
theorem for function fields, due to Geyer and Jarden. The following is a weak consequence
of~\cite[Proposition~13.4]{gejapaper}.

\begin{thm}
\label{gejathm}
Let $c$ be a conjugacy class in $G$. We denote by $P(L/K,c)$ the set of rational
points $x \in \P^1(k)$ outside the branch locus of $C' \rightarrow \P^1_k$
such that \mbox{$c=(x,L/K)$}.
Then one has
\begin{equation}
\label{gejaeq}
\abs{P(L/K,c)} \geq \frac{1}{N}\left( q^s - (N+2g)q^{s/2} - N q^{s/4} - 2(g+N) \right) \text{.}
\end{equation}
\end{thm}

Some preparation is in order before applying Theorem~\ref{gejathm}: to be able to deduce
from it that $P(L/K,c)$ is non-empty as soon as $s$ is chosen large enough, we need to make
sure that the right-hand side of~(\ref{gejaeq}) does grow when $s$ goes to infinity. For instance,
it suffices to establish that $N$ and $g$ are bounded independently of $M$ and $s$.
The integer~$N$ is obviously independent of the choices made: it is equal
to $\prod_{i=1}^n (\deg(f_i)!)$. We shall actually prove that $g$ is also independent
of~$M$ and~$s$.

As $C'$ is geometrically connected over~$k$,
Hurwitz's theorem \cite[IV.2.4]{hartshorne} enables us to express $g$
in terms of the ramification divisor of $C' \rightarrow \P^1_k$.
The finite extension $L/K$ is tamely ramified, since its ramification indices are $\leq 2$
and $p \neq 2$; we can therefore write the ramification divisor in terms of the ramification
indices. We finally obtain the equality
\begin{equation}
\label{eqg}
g-1+N = \frac{N}{2} \sum_{i=1}^n \deg(R_i) \text{,}
\end{equation}
where $R_i \subset \P^1_k$ is now considered as a finite $k$-scheme with its reduced subscheme
structure. Hurwitz's theorem applied to the finite morphism $(C_i)_k \rightarrow \P^1_k$
yields
\begin{equation}
\label{eqr}
\deg(R_i) = 2 g_i - 2 + 2 \deg(f_i) \text{,}
\end{equation}
where $g_i$ denotes the genus of $C_i$. By combining equations~(\ref{eqg}) and~(\ref{eqr}),
we end up with
\begin{equation}
\label{eqgr}
g-1+ N = N \sum_{i=1}^n (g_i - 1 + \deg(f_i)) \text{,}
\end{equation}
hence the claim: $g$ is independent of $M$ and $s$.

We can therefore assume that the right-hand side of~(\ref{gejaeq}) is $\geq 2$, by demanding
that~$s$ be sufficiently large.
The canonical isomorphism $G = G_1 \times \cdots \times G_n =
\symm{H_1 \backslash G_1} \times \cdots \times \symm{H_n \backslash G_n}$ allows us to choose
an element $\sigma \in G$ whose projection
in $G_i$ acts transitively on $H_i \backslash G_i$ for every $i \in \{\tuple{1}{n}\}$.
Let $x_0 \in \P^1(k)$ be the point corresponding to the line in $\P^2_k$ passing
through~$M$ and~$M_0$.
Theorem~\ref{gejathm} now ensures the existence of a rational point $x \in \P^1(k)$
outside $\bigcup_{i=1}^n R_i$, distinct from~$x_0$, and
such that $\sigma=(x,L/K)$. As the image of $(x,L/K)$ in $G_i$ is $(x,K_i/K)$,
it follows from Lemma~\ref{wellknownlemma} and the definition of $\sigma$ that
$\phi_i^{-1}(x)$ is irreducible. Moreover, the $k$-scheme $\phi_i^{-1}(x)$ is étale since
$x \not\in R_i$, and hence it is integral.
That $x \neq x_0$ implies that there exist $a, b \in k$
such that for every $i \in \{\tuple{1}{n}\}$, the scheme
$\Spec(k[t]/(f_i(t,at+b)))$ is an open subscheme of $\phi_i^{-1}(x)$;
as the latter scheme is integral, the polynomials $\tuple{f_1(t,at+b)}{f_n(t,at+b)}$ must
be irreducible.

\section{Application}

The following problem was posed in~\cite{Gao}:

\emph{\textsc{Problem:} Let $f(t,x)\in \Fq[t,x]$ be irreducible and set $g(t)=f(t,at+b)$. Count (or estimate) the number of pairs~$(a,b) \in \Fq \times \Fq$ such that~$g(t)$ is irreducible over~$\Fq$.}

As a partial solution of this problem, we have

\begin{prop}
Let $\Fq$ be a finite field of characteristic $p$ and cardinality $q$.
Let $f \in \Fq[t,x]$ be an irreducible polynomial whose
total degree $d$ satisfies $p\nmid d(d-1)$. Assume that the curve $C \subset \P^2_\Fq$ defined
as the Zariski closure of the affine curve
$$
f(t,x)=0
$$
is smooth and that $q>9(d(d-1)d!+2)^{2}$. Then the polynomial $f(t,at+b) \in \F{q}[t]$
is irreducible for at least $\frac{1}{d!}(q-\frac{d^{4}}{2})(q-3(d(d-1)d!+2)q^{1/2}-d!)$ pairs $(a,b) \in \Fq\times \Fq$.
\end{prop}

\begin{pf}
We unfold the proof of Theorem~\ref{mainthm} for one polynomial and estimate the resulting number of possible pairs $(a,b)$, using only the simplest non-trivial estimates.

The only parts of the proof of Theorem~\ref{mainthm} which depend on the size of the finite field in relation to the degree $d$ are the applications of Proposition~\ref{existsu} and of Theorem~\ref{gejathm}. The case $d=1$ being trivial, we assume $d\geq 2$.

We use the notation introduced in the proof of Theorem~\ref{mainthm}.
As mentioned in the proof of Proposition~\ref{existsu}, the Gauss map~$\phi \colon C \rightarrow
C^\star$ is separable. Proposition 3.5 in~\cite{katz} then implies that  as in characteristic~$0$, the
number~$n$ of lines~$D \subset \P^2_{\Falg}$ whose scheme-theoretic
intersection with $C_{\Falg}$ does not have at most one double point is bounded by the
number of singular points of the dual curve $C^\star_{\Falg}$.
Since $C^\star$ is irreducible of degree $d(d-1)$, we thus obtain $n \leq p_a(C^\star)=\frac{1}{2}
(d(d-1)-1)(d(d-1)-2) \leq \frac{1}{2} d^4-d+1$.
In particular, it follows from our assumption $q>9(d(d-1)d!+2)^{2}$ that $q^2 > n$, hence the existence of an $\F{q}$-rational line $D_0 \subset \P^2_{\F{q}}
\setminus \{M_0\}$
whose scheme-theoretic intersection with~$C$ has at most one double point.

Pairs $(a,b) \in \Fq \times \Fq$ such that $f(t,at+b) \in \Fq[t]$ is irreducible
correspond bijectively to lines~$D \subset \P^2_{\Fq} \setminus \{M_0\}$ such that
$C \cap D$ is integral. Let us denote their number by~$e$, and by $e(M)$
the number of such lines which contain a given $M \in (D_0 \cap U)(\Fq)$.
For such an $M$, we have seen in Theorem~\ref{gejathm} that $e(M)+1$ is greater than or equal
to the right-hand side of~(\ref{gejaeq}) with $s=1$; moreover $N=d!$, and~(\ref{eqgr})
yields
$$
g=1+\frac{N}{2}(d-2)(d+1) \text{,}
$$
whence
$$
e(M) \geq \frac{1}{d!}(q-3(d(d-1)d!+2)q^{1/2}) \text{.}
$$
The result now follows in view of the inequalities
$$
e \geq \sum_{M \in (D_0 \cap U)(\Fq)} \left(e(M)-1\right)
$$
and
$$
\abs{(D_0 \cap U)(\Fq)} \geq q+1-n-d \geq q - \frac{1}{2}d^4 \text{.}
$$
\end{pf}

\section{Appendix}

The following lemma was used several times in the proof of Theorem~\ref{mainthm}.
It is essentially well-known, but we state it here and include a sketch of proof
for lack of an adequate reference. It is only for technical reasons that we
state it in such generality (we need to allow $X'$ to be non-connected in order
to be able to reduce to the case of a strictly Henselian base in the proof below).

\begin{lemma}
\label{wellknownlemma}
Let $u \colon X' \rightarrow X$ and $f \colon X \rightarrow B$ denote
surjective finite flat morphisms
of normal schemes. Assume $B$ is the spectrum of a discrete valuation ring.
Put~\mbox{$f' = f \circ u$}.
Let $G$ be a finite subgroup of $\Aut_B(X')$
such that the generic fibre of $f'$ is a torsor under~$G$. Let $m \in X'$ belong
to the special fibre of~$f'$.
We denote respectively by \mbox{$D_m \subseteq G$} and $I_m \subseteq G$
the decomposition and inertia subgroups associated with~$m$; in other words,
$D_m$ is the stabilizer of $m$ and $I_m$ is the kernel of the natural map
$D_m \rightarrow \Aut(\kappa(m))$.
Let $H = G \cap \Aut_X(X')$.
Then the double quotient $H \backslash G / D_m$ is canonically in bijection with
the special fibre of $f$, and the double quotient $H \backslash G / I_m$ is canonically
in bijection with the geometric special fibre of $f$.
\end{lemma}

\begin{sketchofproof}
Let us first consider the assertion about $H \backslash G / I_m$.
To prove it, one easily checks that $B$ may be assumed to be strictly Henselian,
by using the fact that for any finite field extension $L/k$,
the group $\Aut_k(L)$ acts freely on $\Spec(L \otimes_k \bar{k})$,
where~$\bar{k}$ denotes a separable closure of $k$.
Now the assertion
about~$H \backslash G / I_m$ follows from the one about $H \backslash G/D_m$ since~$D_m=I_m$.

We are thus left with the first part of the lemma. Define a map
$$H \backslash G / D_m \longrightarrow f^{-1}(f'(m))$$ by sending the double class
$H \sigma D_m$ to $u(\sigma(m))$. The key ingredient for
checking that this map is indeed a bijection
is the transitivity of the action of $G$ (resp. $H$)
on the fibres of $f'$ (resp. $u$), and it is a consequence
of~\cite[Ch.~5, §2, Th.~2]{bouralg}.
\end{sketchofproof}

\end{document}